\documentclass{elsart}
\usepackage{natbib}

\usepackage{latexsym}
\usepackage{amsfonts}
\usepackage{amssymb}
\usepackage{amsmath}

\newtheorem{Theorem}{Theorem}[section]

\newtheorem{Problem}{Problem}[section]
\newtheorem{Conjecture}{Conjecture}[section]

\newcommand{\A}{\mathcal{A}}

\newcommand{\CLSYS}{WMIN}
\newcommand{\GWA}{SWR}
\newcommand{\DWA}{WR}
\newcommand{\HDWA}{HDWR}
\newcommand{\HPWA}{HPWR}

\newcounter{codelines}
\newenvironment{sourcecode}{\setcounter{codelines}{1} \centerline{\rule{10cm}{0.5pt}} \medskip \begin{tabular}{ll} }{ \end{tabular} \medskip \centerline{\rule{10cm}{0.5pt}}\setcounter{codelines}{0}}
\newcommand{\codeline}[1]{ \multicolumn{2}{l}{#1} \\}
\newcommand{\codelineone}[1]{\arabic{codelines}: & \addtocounter{codelines}{1}\hskip 1em  \em{#1} \\}
\newcommand{\codelinetwo}[1]{\arabic{codelines}: & \addtocounter{codelines}{1}\hskip 2em  \em{#1} \\}
\newcommand{\codelinethree}[1]{\arabic{codelines}: & \addtocounter{codelines}{1}\hskip 3em  \em{#1} \\}
\newcommand{\codelinefour}[1]{\arabic{codelines}: & \addtocounter{codelines}{1}\hskip 4em  \em{#1} \\}
\newcommand{\reserved}[1]{{\bf \uppercase{#1}}}
\newcommand{\const}[1]{{\bf #1}}
\newcommand{\codecomment}[1]{ /* {\em #1} */}

\begin{document}
\begin{frontmatter}

\title{A Hybrid Search Algorithm for the Whitehead Minimization Problem}

\author{A.D. Myasnikov}
\ead{amyasnik@stevens.edu}
\address{Department of Mathematical Sciences, \\ Stevens Institute of Technology, \\ Hoboken, New Jersey}
\and
\author{R.M Haralick}
\ead{haralick@ptah.gc.cuny.edu}

\address{Department of Computer Science, \\ The Graduate Center of CUNY, \\ New York}

\begin{abstract}
The Whitehead Minimization problem is a problem of finding elements of the minimal length
in the automorphic orbit of a given element of a free group.
The classical algorithm of Whitehead that solves the problem depends exponentially on the group rank.
Moreover, it can be easily shown that exponential blowout occurs when a word of   minimal
length has been reached and, therefore, is inevitable except for some trivial cases.

In this paper we introduce a deterministic Hybrid search algorithm
and its stochastic variation for solving the Whitehead minimization
problem. Both algorithms use search heuristics that allow one to
find a length-reducing automorphism in polynomial time on most
inputs and significantly improve the reduction procedure. The
stochastic version of the algorithm employs a probabilistic system
that decides in polynomial time whether or not a word is
minimal. The stochastic algorithm is very robust.
It has never happened that a non-minimal element has been claimed to be minimal.
\end{abstract}
\begin{keyword}
Automorphism Problem, Free Groups, Heuristic Search.
\end{keyword}

\end{frontmatter}

\section{Introduction}

The Whitehead Minimization problem is a problem of finding elements
of the minimal length in the automorphic orbit of a given element of
a free group. This problem is of great importance in group theory
and topology and continually attracts a great deal  of attention
from the research community.

Starting from the seminal paper of \cite{W}, 
the Whitehead Minimization problem was studied extensively for more
than 70 years (see \cite{LS,Cohen,Lee,Khan,KSS,MS,KKS,K}) and still
the complexity of this problem is unknown.

One of the most important applications of the Whitehead Minimization
problem is that its solution is  part of the solution to the famous
\emph{Automorphism Problem} in free groups introduced by
J.H.C.Whitehead in 1936. Methods used to solve the Whitehead
Minimization problem can be used to decide whether an element is a
part of a generating basis of a free group. The same methods and
their generalizations are used in solving equations over free groups
(see \cite{Razborov}). To practitioners, the Whitehead Minimization
problem could be of interest because of its relation to
non-commutative variations of the public key cryptographic scheme by
\cite{moh}.

All known methods of solving the Whitehead Minimization problem  have exponential dependence
on the rank of a free group. Moreover, the worst case scenario occurs when solving a termination problem (which is to decide
whether or not a given element is minimal) for a minimal element.
Since the goal of the Whitehead Minimization problem is to find a minimal element, the worst case
is inevitable for almost all elements except for elements of a very particular type.
This observation leads us to a conclusion that the known deterministic techniques are
not suitable for groups of large ranks.

  \cite{HMM,HMM2,M,MM}, using methods of pattern recognition and exploratory data analysis,
show that by introducing proper strategies one can construct a length reduction process which is very efficient
on most inputs. Furthermore,  in these papers we formulate several conjectures (see \ref{sec:WMIN}) regarding the various properties of the problem.

In this paper we present a new algorithm  for solving the Whitehead Minimization
problem. It is a \emph{hybrid} algorithm in a sense that it employs several
 stochastic, as well as deterministic, procedures  based on the conjectures stated by \cite{HMM}.

We combine a stochastic search algorithm  and heuristic search procedures
(both described in Section \ref{sec:heur}) with the probabilistic
classification system ``recognizing" minimal elements (see Section \ref{sec:gaus})
to construct a Hybrid Deterministic Whitehead
Reduction (\HDWA) algorithm  solving the Length Reduction Problem in
a polynomial number of steps (in terms of  group rank) on most input words from a free group.
The resulting algorithm is  deterministic and  still requires an exponential number
of steps to prove that a word is minimal.

We present a fast probabilistic algorithm \HPWA
which is a slight modification of \HDWA. Algorithm  \HPWA\ is very robust and extremely fast
on most input words, including words in free groups of large ranks. Although we do not have a
formal proof of the correctness of \HPWA, in all the experiments that we have performed it has never
happened that the algorithm  has produced an incorrect output.

The algorithms \HDWA\ and its probabilistic version \HPWA\ are described in Section \ref{sec:alg_descr}.
We  give experimental results evaluating the performance
of these two algorithms in Section \ref{sec:HDWA_exps}. Comparison with the standard
deterministic procedure is also presented in Section \ref{sec:HDWA_exps}.

\section{The Whitehead minimization problem}
\label{sec:WMIN}

\label{sub-sec:WminProb} \label{sec:WDA}

Let $X = \{x_1, \ldots, x_n\}$ be a finite alphabet; $X^{-1} =
\{x^{-1} \mid x \in X\}$ be the set of formal inverses of letters
from $X$ and $X^{\pm 1} = X \cup X^{-1}$. A word $w = y_1 \ldots y_m
$ in the alphabet $X^{\pm 1}$ is called {\em reduced} if $y_i \neq
y_{i+1}^{-1}$ for $i = 1, \ldots, m-1$ (here we assume that
$(x^{-1})^{-1} = x$).    Applying reduction rules $xx^{-1}
\rightarrow \varepsilon, x^{-1}x \rightarrow \varepsilon$ (where
$\varepsilon$ is the empty word), one can reduce each word $w$ in
the alphabet $X^{\pm 1}$ to a reduced word $\overline {w}$. The word
$\overline {w}$ is uniquely defined and does not depend on a
particular sequence of reductions. Denote by $F = F(X)$ the set of
reduced words over $X^{\pm 1}$. The set $F$ forms  a group with
respect to the multiplication $u \cdot v = \overline{uv}$,
called   a {\em free } group with {\em basis} $X$. The cardinality
$|X|$  is called the {\em rank}
 of $F(X)$. We write $F_n$ instead of $F$ to indicate that
 the rank of $F$ is equal to $n$.

A bijection $\phi: F \rightarrow F$ is called an {\em automorphism}
of $F$ if $\phi(uv) = \phi(u)\phi(v)$ for every $u, v \in F$. The
set  $Aut(F)$ of all automorphisms of $F$ forms a group  with
respect to  the composition of maps. Every automorphism $\phi \in
Aut(F)$ is completely determined by the images $\phi(x)$ of elements
$x \in X$. Sometimes it is more convenient to  use non-functional notation $w\phi$ to
denote the action of automorphism $\phi$ on $w$.

The following two subsets of $Aut(F)$ play an important part in both group
theory and topology.
 An automorphism  $t \in Aut(F)$  is called a
 {\em Nielsen automorphism} if  for some $x \in X$ $t$ fixes
 all elements $y \in X, y \neq x$ and maps $x$ to one of the
 elements $x^{-1}$, $y^{\pm 1}x$, $xy^{\pm 1}$. Note that
 automorphisms that map $x$ to $x^{-1}$, leaving everything else unchanged, cannot cause alterations of
 the word length. Such automorphisms will be called length invariant
 automorphisms. By $N(X)$ we
 denote the set of all Nielsen automorphisms of $F$ except the length-invariant ones.

 A non-trivial automorphism  $t \in Aut(F)$  is called a
 {\em Whitehead automorphism} if it has one of the following types:
 \begin{list}{}{}
  \item 1) $t$ permutes the elements of $X^{\pm 1}$;
  \item 2)   $t$ fixes a given element $a \in X^{\pm 1}$ and maps each element $x
\in X^{\pm 1}, x \neq a^{\pm 1}$ to one of the elements $x$, $xa$,
$a^{-1} x$, or $a^{-1} x a$.
\end{list}
It is easy to see that automorphisms of the first type are length-invariant. By $W(X)$ we denote the set of
Whitehead's automorphisms of the second type.
Obviously, every Nielsen automorphism
is also a Whitehead automorphism.

Observe that
  $$|N(X)| = 4n(n-1), \ \ \ |W(X)| = 2n4^{(n-1)} - 2n$$
where $n = |X|$ is the rank of $F$.

It is known (see \cite{LS}) that every automorphism from
$Aut(F)$ is a product of finitely many Nielsen (hence Whitehead)
automorphisms.

 The automorphic orbit
$Orb(w)$ of  a  word $w \in F$ is the set of all automorphic
images of $w$ in $F$:
\[Orb(w) = \{ v \in F  \mid  \exists \varphi \in Aut(F) \mathrm{\; such \;that}\; w\varphi =
v\}.\] A word $w \in F$ is called {\em minimal} (or {\em
automorphically minimal}) if $|w| \leq |w\varphi|$ for any
$\varphi \in Aut(F)$. By $w_{min}$ we denote a word of minimal
length in $Orb(w)$. Notice that since there may be several elements of the minimal
length in same orbit, $w_{min}$ is not unique in general.

\begin{Problem}[Minimization Problem (MP)]
\label{RP} For a word $u \in F$ find an automorphism $\varphi \in
Aut(F)$ such that $u\varphi = u_{min}$.
\end{Problem}

 \cite{W}   proved the following result which gives
 a solution to the minimization problem.

\begin{Theorem}[Whitehead]
\label{thm:T1} Let $u \in F_n(X)$. If
$|u|>|u_{min}|$, then there exists $t \in W(X)$ such that
\[
|u|>|ut|.
\]
\end{Theorem}

An automorphism $\phi \in Aut(F)$  is called {a length-reducing}
automorphism for a given  word $u \in F$ if $|u\phi| < |u|$. The
theorem  above claims that the finite set $W(X)$ contains a
length-reducing automorphism for every non-minimal word $u \in F$.
This allows one to design a simple search algorithm for (MP).

 Let $u \in F$. For each $t \in W(X)$ compute the length
of the word  $ut$ until $|u|
> |ut|$, then put $t_1 = t, u_1 = ut_1$. Otherwise stop and output
$u_{min} = u$. This procedure  is called the {\em Whitehead
Length Reduction}  routine (WLR). Now  Whitehead Reduction (\DWA)
algorithm  proceeds as follows. Repeat WLR on $u$, and then
on the resulting $u_1$, and so on, until at some  step $k$ WRL gives an output
$u_{min}$. Then $ut_1 \ldots t_{k-1} = u_{min}$, so $\phi = t_1
\ldots t_{k-1}$ is the required automorphism.

 Notice,  that the iteration procedure \DWA\   simulates the classical
 greedy  descent method ($t_1$ is a successful  direction from $u$;
 $t_2$ is a successful direction from $u_1$;  etc.) Theorem \ref{thm:T1}
guarantees that the greedy approach will always converge to a
global minimum.

Clearly, there could be at most $|u|$ repetitions of WLR on an
input $u \in F$
\[|u| > |ut_1|>...>|ut_1 ... t_l| = u_{min}, \ \ \ l \leq |u|.\]
 Hence the worst case complexity of the Whitehead's algorithm is bounded from above by
\[
cA_n|u|^2,
\]
where $A_n = 2n4^{(n-1)} - 2n$ is the number of Whitehead
automorphisms in $W(X)$ and the constant $c$ is a stretching factor
by which the length of a word increases  after a Whitehead
automorphism is applied (ignoring the low level implementation
details.) One letter can be mapped into a word of length of at most 3, so $c$ is
bounded by 3 and does not depend on the rank of a group or the
word's length. Since $A_n$ depends exponentially on the rank of a
free group, in the worst case scenario the algorithm seems to be
impractical for free groups with large ranks. One can try to improve
on the number of steps which it takes to find a length-reducing
automorphism for a given non-minimal element from $F$. In this
context, the  question of interest is the complexity of the
following

\begin{Problem}[Length Reduction Problem(LRP)]
For a given non-minimal element $u \in F$ find a length-reducing
automorphism.
\end{Problem}

We refer to \cite{HMM,MM} for a general discussion of this problem.
 \cite{HMM}   offers some empirical
evidence that by using  smart strategies in selecting Whitehead
automorphisms $t \in W(X)$  one can dramatically improve the
average complexity of WLR in terms of the rank of a group.
Some of the experimental results were formulated as the following conjectures:
\begin{Conjecture}[\cite{HMM}]
\label{conj:NielsenFirst}
 Let $U_k$ be the set of all non-minimal elements in $F$ of length $k$
and $NU_k \subset U_k$ the subset of elements  which have Nielsen
length-reducing automorphisms. Then
 $$
\lim_{k \rightarrow \infty} \frac{|NU_k|}{|U_k|}  = 1.
$$
\end{Conjecture}

\begin{Conjecture}[\cite{HMM}]
\label{conj:hypersurf}
The feature vectors of weights of the Whitehead Graphs of elements
from $F$ are separated into bounded  regions in the corresponding
 space. Each such  region can be bounded by a hypersurface and
corresponds to a particular Nielsen automorphism  in a sense that
all elements in the
 corresponding class can be reduced by that automorphism.
 \end{Conjecture}

Arguably the conjectures above are not intuitive and most likely would have
been difficult to arrive at without observations obtained using
computer experiments. At this point we would like to mention a
new development in this area which was not available during the
submission of this paper. \cite{K} has recently posted a preprint,
giving a mathematical proof of the Conjecture
\ref{conj:hypersurf}. To our best knowledge this is the first time
 non-trivial conjectures were obtained using
statistical and exploratory data analysis techniques.

Unfortunately, one can easily see that the worst case behavior of the algorithm \DWA\
occurs when a word of minimal length has been reached.
Except for some trivial cases (when a minimal word is a generator, for example)
all Whitehead automorphisms need to be applied to a minimal word  before we can conclude that
it is, indeed,  minimal. It seems that no algorithm is known to  avoid  time-consuming
computation in this case. We would like to emphasize the importance of this
fact by formulating it as a separate problem:
\begin{Problem}[Minimal Word Classification  Problem(MWCP)]
For a given $u \in F(X)$, decide whether $u$ is minimal or not.
\end{Problem}
We discuss this problem in the previous papers. \cite{HMM2}  gives a probabilistic solution which
is based on regression models. \cite{M} used the so-called support vector machines to improve the
performance in free groups of large ranks. In this paper
we introduce a new,  significantly more efficient probabilistic system based
on the empirical distribution of minimal elements in the corresponding feature space (see Section \ref{sec:gaus}).

\section{Description of the Hybrid algorithms}
\label{sec:alg_descr}
\subsection{Heuristics for the Length Reduction Problem}
\label{sec:heur}

We have addressed this problem in the preceding papers. Our first
approach was to develop a simple Stochastic Whitehead Reduction (\GWA) algorithm to solve LRP. It is
implemented as a combination of a greedy descent procedure with  genetic
search techniques.

Define the search space $\mathcal{S}$  as the set of all finite
sequences \[\mu  = <t_1, \ldots, t_s>\] of Whitehead automorphisms $t_i \in W(X)$.
For such $\mu$ and a word $u \in F$  define $u\mu = ut_1 \ldots t_s$.

The solution to LRP is any sequence $\mu^* \in \mathcal{S}$ such that
\[|u\mu^*| < |u|.\]
Among all such solutions we prefer the ones that give maximal
length reduction of the image. In \GWA, we define the criterion
function which evaluates a solution $\mu$ as
\[\mathcal{F}(\mu) = |u\mu|.\]
The details on the implementation and evaluation of \GWA\ can be found
in \cite{MM}.

To our great surprise, this naive stochastic algorithm significantly outperformed the standard
algorithm, especially in free groups of large ranks. For example, there were very few runs
of \DWA\  for words $w \in F_{10}$ with $|w| > 100$ that finished within an hour and there
were no such runs for $|w| > 200$. Nevertheless, the stochastic algorithm still was able to
find minimal words in a matter of seconds. What seemed to be more important is that the stochastic
algorithm did not show exponential dependence on the group's rank.

We  strongly believe that if a stochastic algorithm performs very well, then there must
be a purely mathematical reason behind this phenomenon which can be uncovered by a proper
statistical analysis. Following this philosophy, we performed an analysis of successful  solutions
produced by \GWA. The results helped us to define a number of search heuristics described by \cite{HMM}.
Below we give a brief description of these heuristics.

First, we observed that among all Whitehead automorphisms in the successful solutions,
Nielsen automorphisms statistically had a greater chance to occur. Further experiments showed
that more than 99\% of non-minimal elements can be reduced by one of the Nielsen automorphisms.
Our first heuristic is based on this observation and simply
suggests trying Nielsen automorphisms first in the routine WLR,
i.e., in this case we assume that in the fixed listing of
automorphisms of $W(X)$, the automorphisms from $N(X)$ come first.
We refer to this heuristic as  \emph{Nielsen First}.
Note, that the Nielsen First heuristic is very general and does not use information about the input word itself.
We showed that one can significantly improve the performance of the search procedure by incorporating
heuristics that use some knowledge about the input.

Let $u \in F(X)$. The undirected {\em Labeled Whitehead graph} $W(u)= (V,E(u))$ of the
word $u$ is a complete undirected graph, where the set of vertices $V$ is equal to the set  $X^{\pm 1}$.
Every edge $e = (x,y)$, $x \neq y$  of the Whitehead graph is assigned a weight $\omega_e = n_e / |u|$, where $n_e$ is the number of times
subwords $x y^{-1}$ or $y x^{-1}$ occur in $u$. Note that $\omega_e = 0$ if the subwords corresponding to the edge $e$ do not occur in $v$.
Now, for a given word $u \in F$ define a special vector representation (called a \emph{feature vector}) $f(u) \in \mathbb{R}^{|E(u)|}$ such that
\[f(u) = < \omega_{e_1}, \ldots, \omega_{e_{|E(u)|} } >.\]
The edges $e_i$ are assumed to be taken in some fixed order. Since the Whitehead graph is complete,
the number of edges and, therefore, the size of feature vectors is $3n^2 -n$ for all
elements in a free group $F_n$.
The set of all feature vectors is usually called \emph{a feature space} and is denoted by $\mathcal{F}$.

Experiments show that there is a correlation between the location of the feature vectors in the corresponding
space and the length-reducing Nielsen automorphisms.

Let $t \in N(X)$ be a Nielsen automorphism. Define the set
$$O_t = \{ w \mid  r \in N(X)\ \mathrm{and} \   |wr| < |w| \Longleftrightarrow r = t \}$$
as a set of all  elements   that can be
reduced only by $t$ and no other Nielsen automorphism. We also define a set $B_{m,t} \subset O_t$:
$$B_{m,t} = \{ w \mid  w \in O_t, |w| \leq m\},$$
which is a finite set of elements from $O_t$ with the length of at most $m$.
For a large $m$ we define
$$\lambda_t = \frac{1}{|B_{m,t}|} \sum_{w \in B_{m,t}} f(w)$$
as an estimate of the {\em mean feature vector} of the elements in $O_t$.


Now, let
$$d(w,t) = ||f(w) - \lambda_{t}||$$
be the distance (in this case, Euclidean distance) between the feature vector of
a given word $w$ and the estimate $\lambda_t$ of the mean feature vector corresponding to the Nielsen automorphism $t$.

\cite{HMM}  experimentally  show that in about 99\% of the time a randomly generated non-minimal element $w$
can be reduced by a Nielsen automorphism $t^*$ such that
\[
d(w,t^*) = \min\{ d(w,t) \mid t \in N(X) \}.
\]

Now we define the second heuristic, which is called the  \emph{Centroid} heuristic. For a given
word $u$ compute distances $d(u,t)$ for all $t \in N(X)$ and sort them in  the increasing order:
$$d(u,t_1) \leq d(u,t_2) \leq \ldots \leq d(u,t_k).$$ Apply  automorphisms $t_{1}, \ldots, t_{k}$ sequentially until
a length-reducing Nielsen automorphism (if any) is found.

Now let $e = (x, y^{-1})$, $x, y^{-1} \in V$  be an edge in the Whitehead graph of a word $w$ such that  $x \neq y$.
By construction of the Whitehead graph, $e$ corresponds to subwords $s_e = (x y)^{\pm 1}$.

There are  only two  Nielsen automorphisms that reduce lengths of subwords in $s_e$:
$$
\psi^x_e: x \rightarrow x y^{-1},\ z \rightarrow z \ \forall  z \neq x
$$
and
$$
\psi^y_e: y \rightarrow x^{-1} y,\ z \rightarrow z \ \forall  z \neq y.
$$
Denote $\psi_e = \{ \psi^x_e, \psi^y_e \}$. We call  automorphisms $\psi_e$ the  {\em length reducing with respect} to
the edge  $e$.

We can order Nielsen automorphisms $\psi_{e_i} \subset N(X)$:
\begin{equation}
\label{eq:max_weight_order}
<\psi_{e_1}, \psi_{e_2}, \ldots, \psi_{e_k}>
\end{equation}
such that the corresponding edges $e_1, \ldots , e_{k}$ are chosen
according to  the decreasing order of the values of their weights
\[
\omega(e_1) \geq \omega(e_2) \geq \ldots \geq \omega(e_{k}).
\]

In the third heuristic we apply Nielsen automorphisms in the
order given by (\ref{eq:max_weight_order}). This heuristic is
called the Maximal Edge heuristic. In the paper (\cite{HMM}) we present  empirical evidence
that most non-minimal elements can be reduced by one of the automorphisms in
$\psi_{e_1}$, given that $\omega(e_1)$ is maximal.

To get a better understanding of how effective these methods are, we estimate the 99th percentile of the number
of Nielsen automorphisms required to reduce a  non-minimal word from a given test set (see Table \ref{tab:99PC}).

\begin{table}
\begin{center}
\begin{tabular}{|l|c|c|c|c|}
\hline
  & $F_3$ & $F_4$ & $F_5$ \\
\hline
$|N|$ &  24 & 48  & 80  \\
\hline
$S_1$ &  1  &  1 &   1 \\
\hline
$S_{10}$ & 1  & 2  &  3  \\
\hline
\multicolumn{4}{c}{Centroid} \\

\end{tabular} \ \
\begin{tabular}{|l|c|c|c|c|}
\hline
  & $F_3$ & $F_4$ & $F_5$ \\
\hline
$|N|$ &  24 & 48  &  80  \\
\hline
$S_1$ &  2  & 3  &   4 \\
\hline
$S_{10}$ & 6  & 7  &  6  \\
 \hline
\multicolumn{4}{c}{Maximal Edge} \\
\end{tabular}\ \
\begin{tabular}{|l|c|c|c|c|}
\hline
  & $F_3$ & $F_4$ & $F_5$ \\
\hline
$|N|$ & 24 & 48 & 80 \\
\hline
$S_1$ & 23  & 45 &  70 \\
\hline
$S_{10}$ & 21 & 37 &  50 \\
 \hline
\multicolumn{4}{c}{Nielsen First} \\
\end{tabular}
\end{center}
\caption{ 99th percentile of the number of Nielsen automorphisms computed for different heuristics applied to test sets $S_1$ and $S_{10}$ in free groups $F_3$, $F_4$ and $F_5$. }
\label{tab:99PC}
\end{table}

We can see that the Centroid heuristic is the most effective and is able to predict a length reducing automorphism with very high
accuracy. The Maximal Edge heuristic also uses very few automorphisms
where the  Nielsen First heuristic must apply at least 60\% of
Nielsen automorphisms in the best case.

The Nielsen First heuristic does not require
any additional computations and, therefore, its computational complexity is of order $\mathrm{O}(1)$.
The Maximal Edge heuristic requires $\mathrm{O}(|w| + n^2)$ steps. The Centroid heuristic  requires
$\mathrm{O}(|w| + n^4)$ elementary steps, where $n$ is the rank of a free group.
Therefore, one becomes aware of a tradeoff between the length of the input word and the rank of
a group. Since the Centroid and Maximal Edge heuristics are more accurate they become more attractive when the length of the input word increases because
fewer superfluous automorphisms will be applied to the input word.

In Section \ref{sec:HDWA} we show how the stochastic algorithm SWR can be combined
with Centroid and Maximal Edge heuristics to improve the solution of the Whitehead
Minimization problem.

\subsection{Probabilistic System for Classification of Minimal Words}
\label{sec:gaus}
We have already mentioned that the worst case of the standard Whitehead algorithm
applied to solve LRP occurs when the word is already minimal. Note
that  exponential blowout is
inevitable in the WR algorithm, unless minimal words are of  a very special type.  Being able to
solve the Minimal Word Classification Problem efficiently is crucial for an efficient solution
to the Whitehead Minimization Problem. In \cite{HMM2} and \cite{M} we describe several stochastic classification
systems (classifiers) based on pattern recognition techniques such as regression and support vector machines.
These classifiers are able to decide whether a given word is minimal   in polynomial time (with respect to
group rank) with a very small error of misclassification.

Conclusions in \cite{M} suggest that one of the classes (minimal or non-minimal) of elements
could be located in a compact region in the feature space $\mathcal{F}$ and can be bounded by a
hypersurface.

To support this conjecture we perform the following experiment. Assume that the feature
vectors of minimal elements follow the multivariate normal distribution $\mathcal{N}(\mu,\Sigma)$ with the
mean $\mu$ and the covariance $\Sigma$. We estimate $\mu$ and $\Sigma$ from a set of randomly generated
minimal elements. Experiments show that more than 97\% of minimal elements lie inside the hyperellipse,
corresponding to the 99.9\% confidence interval for $\mu$. Moreover, no non-minimal elements
fall inside that region. This is a very strong indication that the feature vectors of minimal elements indeed lie compactly
in $\mathcal{F}$.

Using this result we construct a new probabilistic system \CLSYS\ to solve the Minimal Word Classification Problem.
We decide that a given word $u$ is minimal if its feature vector $f(u)$ falls inside the corresponding hyperellipsoid and
we decide that $u$ is otherwise non-minimal. To be more precise, let $\mu$ and $\Sigma$ be, respectively, the
mean  and the covariance matrix of feature vectors of minimal elements. Using the
so-called Mahalanobis distance we define the decision rule:
\[
decide(u) = \left \{
\begin{array}{ll}
\mathrm{minimal}, & \mathrm{if} \ (x-\mu)^T\Sigma^{-1}(x-\mu) < \rho; \\
\mathrm{non-minimal}, & \mathrm{otherwise}
\end{array}
\right.
\]
where $(x-\mu)^T$ is the transpose of a column vector $x-\mu$.
One way of estimating the threshold $\rho$ was indicated above, where it was taken to
correspond to the 99.9\% confidence interval of $\mu$, given that feature vectors follow
multivariate normal distribution. However, in this case the error of misclassifying
minimal  elements is unacceptably large (greater then 5\%). This indicates that the feature
vectors actually are not normally distributed.

A practical way to estimate $\rho$ is to estimate the distribution of distances from feature vectors
of minimal elements to their mean. Then we take $\rho$ such that $100(1-\alpha)$ percent of  minimal elements
have distances less than $\rho$ for a given $\alpha$. Note that $\alpha$ corresponds to a confidence
level in a non-parametric hypothesis testing.

To compute Mahalanobis  distance we need to obtain $\mu$ and $\Sigma$. One way is to estimate them from a
set of randomly generated minimal elements. This process is usually called ``training" the classifier.
Unfortunately, to generate the sample of minimal elements we require to solve the length reduction problem
which, as we have argued, is hard in groups of large ranks. Below we suggest a more efficient training procedure.

\cite{KSS} show that a random cyclically reduced element in a free group is minimal with
asymptotic probability 1. It is also easily shown that any minimal element is already cyclically reduced.
Following these two facts, we suggest  estimating $\mu$ and $\Sigma$ from a set of randomly generated elements
of a whole free group. This can be implemented very efficiently even in groups with large ranks (see \cite{MM} for
more details on generating random elements in a free group).

We would like to mention here the two kinds of errors that may occur when solving the Minimal Word Classification Problem.
The first is the error  of classifying  a non-minimal element as minimal. It is called the \emph{false positive} error.
The second is the  error of classifying  a minimal  element  as a non-minimal element. It is called the \emph{false negative} error.
The rate of the false positive error in all our experiments was zero. This property of the classifier is very important
for a successful implementation of the probabilistic version of the hybrid algorithm. We will return to this discussion
when we describe \HPWA\ in Section \ref{sec:HPWA}.

\subsection{\HDWA}
\label{sec:HDWA}

The  deterministic hybrid procedure \HDWA\
is given in Figure \ref{fig:HDWA}. The algorithm contains two
major parts. The first part consists of a number of so-called \emph{fast checks} --
linear  or polynomial procedures that can solve the length reduction problem
on  some inputs. In fact, the fast checks used in \HDWA\ are expected to reduce most
non-minimal elements in a free group. The problem is that there are
non-minimal  words which cannot be reduced by fast procedures. Using the fast checks alone, one
cannot  decide whether an input word is minimal or not. We need to provide a termination  condition
of the algorithm. This task is solved by the second part of the algorithm which is a version
of the standard deterministic algorithm \DWA. Note that in most cases, the computationally ineffective
procedure \DWA\ is expected to be executed only on minimal elements.

\begin{figure}[t]
\begin{center}
\footnotesize

\begin{sourcecode}
\codeline{\reserved{Deterministic Hybrid Algorithm} }
\codelineone{\reserved{set} the current word $w_c$; $reduced$ = \const{true};}
\codelineone{\reserved{while} $reduced$ \reserved{begin}}
\codelinetwo{$reduced$ = \const{false};}
\codelinetwo{\codecomment{Begin fast checks}}
\codelinetwo{\reserved{if} $w_c$ \reserved{NOT} classified as  minimal \reserved{begin}}
\codelinethree{\reserved{IF} $\psi_{e_{max} }$ (Maximal Edge) reduces $w_c$}
\codelinefour{$w_c$ = reduced word; $reduced$ = \const{true};}
\codelinethree{\reserved{else if} Centroid reduces $w_c$}
\codelinefour{$w_c$ = reduced word; $reduced$ = \const{true};}
\codelinethree{\reserved{else if} Stochastic Algorithm reduces $w_c$}
\codelinefour{$w_c$ = reduced word;  $reduced$ = \const{true};}
\codelinetwo{\reserved{END IF}}
\codelinetwo{\codecomment{End fast checks}}
\codelinetwo{\reserved{IF NOT}  $reduced$ \reserved{AND} $W(X)$ reduce $w_c$}
\codelinefour{$w_c$ = reduced word; $reduced$ =  \const{true};}
\codelineone{\reserved{END WHILE} }
\codelineone{\reserved{STOP};}
\end{sourcecode}
\end{center}
\caption{Algorithm \HDWA.}
\label{fig:HDWA}
\end{figure}
%
%

\HDWA\ is an iterative procedure. On each iteration, the length reduction problem is
solved for the word $w_c$, which is an  automorphic image of the minimal length
of the input word $w$ found so far.
The algorithm terminates when there are no reductions possible and the current word $w_c$
is returned as a minimal word $w_{min}$.

The first step in the algorithm  is the classification
procedure \CLSYS\ (line 5) which decides whether a current word is minimal or not.
Even though reduction procedures used as fast checks do not require significant computational resources,
this step   helps avoiding superfluous computations by distinguishing minimal elements on the first stage
of the algorithm.

Fast reduction procedures  are based on the search heuristics described
in Section \ref{sec:heur}. We use the Maximal Edge heuristic as the first fast check because
 it requires the least number of steps
($\mathrm{O}(n^2 + |w|)$) when compared to other methods. Moreover, $n^2$ part appears when we construct
the  Whitehead graph which is required  for all heuristics. Note that
more than 90\% of non-minimal elements are expected to be reduced using one of the two
automorphisms corresponding to the maximal weight edge of $WG(w_c)$.

Let $e_{max}(w_c)$ be the maximal edge in $WG(w_c)$ and
$\psi_{e_{max}}$ be the set of two length-reducing automorphisms with respect to
the edge $e_{max}(w_c)$. On line 6 of the algorithm \HDWA\ we apply automorphisms
$\psi_{e_{max}}$ to $w_c$.
The maximal number of steps required to perform the fast check is
\[
\mathrm{O}(n^2 + |w_c|).
\]
These steps include the construction of the feature vector
$f_{WG}(w_c)$  and the application of the automorphisms $\psi_{e_{max}}$.
Following the observations from \cite{HMM}
we expect most  non-minimal elements to be reduced at line 6. If the word $w_c$
has been reduced by one of the two  automorphisms, say $\psi'_{e_{max}}$, we
substitute $w_c$ with $\psi'_{e_{max}}(w_c)$ and start a new iteration.
If the Maximal Edge check fails, we continue the reduction process utilizing the next step.

As the next fast check we use the Centroid heuristic.
Let $\psi_{e_{max}}$ be the set of two automorphisms applied at line 6 and $N(X)$ be the
Nielsen. We order automorphisms
\begin{equation}
\label{eq:HDWA_centr_order}
<\varphi_1, \ldots, \varphi_k>, \ \varphi_i \in N_n - \psi_{e_{max}},
\end{equation}
such that
\[
d(w_c,\varphi_1) \leq d(w_c,\varphi_2) \leq \ldots \leq d(w_c,\varphi_k).
\]

We apply automorphisms $\varphi_1, \ldots, \varphi_k$ in the order given by (\ref{eq:HDWA_centr_order})
to the word $w_c$. If one of the automorphisms has reduced the length of $w_c$, we stop
and start a new iteration   with the new, reduced word $w_c$. The maximal
number of steps required to execute this fast check is
\[
\mathrm{O}(n^4 + n^2|w_c|).
\]

By line 10 we already know that the word $w_c$ does not have Nielsen length-reducing automorphisms.
The suggested strategy at this point is to try to reduce $w_c$ by executing  \GWA\ for a predefined number of generations.
If \GWA\ fails to find a length-reducing automorphism, we continue with the algorithm \DWA.
A very conservative bound for the expected maximal number of generations of the stochastic algorithm
was given in \cite{MM}.
Note that since algorithm \GWA\ performs better than \DWA\ only in groups with relatively large ranks (greater then 5),
it might happen that performance improves if we omit step 10 when the rank of a free group is small.

The maximal time complexity to find a length-reducing automorphism for $w_c$ using \HDWA\ is still
\[
O(2^n |w_c|).
\]
However, following the   discussion in Section \ref{sec:heur},  we
expect the length reduction process to be extremely efficient for most non-minimal words.
Unfortunately, as previously stated, the worst case behavior of the algorithm occurs when the current word $w_c$
becomes minimal and, therefore, it is inevitable except for some trivial cases.
In the next section we introduce a probabilistic algorithm that addresses this problem.

\subsection{\HPWA}
\label{sec:HPWA}

Words of both types, minimal and non-minimal,  may cause an exponential blowout in the
algorithm \DWA. However, non-minimal words  do not seem to be a major problem since we have shown
that most of them can be reduced by one of the Nielsen automorphisms.

The bottleneck in solving the length reduction problem
occurs in  the lack of a fast algorithm to decide whether a word is minimal or not.
In fact, the only known deterministic solution is the algorithm \DWA\ itself.
Recall that the worst case of the algorithm occurs when the input word is
already minimal.  In this case all of the Whitehead automorphisms have to be applied
to the word before the decision that the word is minimal can be made.

In this section we introduce a Hybrid Probabilistic Whitehead Reduction algorithm \HPWA\ for
solving Whitehead's Minimization problem. In \HPWA\ the decision on whether or not
a word is minimal is made using a probabilistic classification system \CLSYS.
This allows one to avoid the exponential blowout for the cost of a possibility of a very small classification error.

\begin{figure}[t]
\begin{center}
\footnotesize
\begin{sourcecode}
\codeline{\reserved{Stochastic Hybrid Algorithm} }
\codelineone{\reserved{set} the current word $w_c$; $reduced$ = \const{true};}
\codelineone{\reserved{while} $reduced$ \reserved{begin}}
\codelinetwo{$reduced$ = \const{false};}
\codelinetwo{\reserved{IF} $\psi_{e_{max} }$ (Maximal Edge) reduces $w_c$}
\codelinethree{$w_c$ = reduced word; $reduced$ = \const{true};}
\codelinetwo{\reserved{else if} Centroid reduces $w_c$}
\codelinethree{$w_c$ = reduced word; $reduced$ = \const{true};}
\codelinetwo{\reserved{else if} $w_c$ classified as minimal }
\codelinethree{\reserved{STOP};}
\codelinetwo{\reserved{else if} Stochastic Algorithm reduces $w_c$}
\codelinethree{$w_c$ = reduced word;  $reduced$ = \const{true};}
\codelineone{\reserved{END WHILE} }
\codelineone{\reserved{STOP};}
\end{sourcecode}
\end{center}
\caption{Algorithm \HPWA.}
\label{fig:HPWA}
\end{figure}
%
%

We construct \HPWA\ from \HDWA\ first by removing the last step (line 14) from the algorithm (see Figure \ref{fig:HPWA}).
Note the increased role the stochastic algorithm \GWA plays. This is the only method in \HPWA\
which is capable of reducing non-minimal elements that do not have Nielsen length-reducing
automorphisms.

Secondly, we move the classification step  behind  the fast reduction procedures. To explain
this modification we would like to return to
the discussion of the roles played by the two types of classification errors of the classifier \CLSYS\ in the
view of its new application.

Recall that the two  errors are: the false positive error (classifying  a non-minimal element as minimal) and
 the false negative error (classifying  a minimal  element  as a non-minimal element).
Now observe that once the classifier \CLSYS\  decides that the word $w_c$ is minimal,
algorithm \HPWA\  terminates and returns $w_c$ as the result. There is no backtracking or additional
checking performed after the decision is made. This means that if a non-minimal word $w_c$ is classified as minimal
the algorithm will produce an incorrect result.
On the contrary, when a minimal word is misclassified as non-minimal, the cost
of such  error is the number of  extra computational steps performed by the algorithm in order to reduce
a non-reducible word. What is important is that the algorithm still produces a correct result.

Let $\epsilon$ be the probability of committing  the false positive error by \CLSYS. Now assume that
during the reduction process classifier \CLSYS\ was called $k$ times to decide whether an element
$w_c$ is minimal or not. The probability that the algorithm terminates with a correct answer
is $(1 - \epsilon)^k$. This shows that the probability of giving an incorrect answer grows rapidly
with the number of times the minimality decision is made.

Note that most of the reductions are expected to be done by fast check procedures.
Moving the classification step behind the fast checks allows us to reduce the error of producing an incorrect
answer while still maintaining a small computational cost on average.

The arguments above show that the false positive error of the classifier has crucial importance.
It is necessary to keep the rate of the false positive
error as minimal as possible in order for the algorithm to perform correctly.

It has been noted in Section \ref{sec:gaus} that  the error of misclassifying non-minimal elements
was zero in all our experiments and, therefore, we expect it to be very small in all instances.

\section{Evaluation}
\label{sec:HDWA_exps}

In this section we  evaluate the  algorithms \HDWA\ and \HPWA\ and compare their performance
to the performance of the algorithm \DWA.

We evaluate these algorithms on the following test sets of randomly generated cyclically reduced non-minimal elements:
 \begin{enumerate}
 \item[$S_1$:] contains minimal and non-minimal elements in equal proportions. Non-minimal elements
are obtained with one Whitehead automorphism.

\item[$S_P$:] set of pseudo-randomly generated   \emph{primitive elements} in $F$. Recall, that
$w \in F(X)$ is primitive if and only if there exists an automorphism $\alpha \in Aut(F)$ such that  $\alpha(w) \in X^{\pm}$.

\item[$S_{10}$:]  generated similarly to $S_1$, but up to 10 automorphisms are used to
generate non-minimal elements.
 \end{enumerate}
Some characteristics of the sets in free groups $F_3$, $F_4$ and $F_5$ are given in  Table \ref{tab:heur_data_non-min}.

%
\begin{table}[t]
\begin{center}
\begin{tabular}{|c|c|c|c|c|c|}
\hline
 Dataset  & Size & Min. length& Avg. length &Max. length & Std. deviation \\
\hline
$S_1$ & 10143 & 3 & 605.8 & 1306 &  359.3\\
\hline
 $S_{10}$ & 2535 &  3 & 1507.9 & 13381 & 1527.9 \\
\hline
 $S_P$ & 5645  & 3 & 1422.1 & 143020 & 5379.0\\
\hline
\end{tabular}

\medskip

a) $F_3$;

\bigskip

\begin{tabular}{|c|c|c|c|c|c|}
\hline
 Dataset  & Size & Min. length& Avg. length &Max. length & Std. deviation\\
\hline
 $S_1$  &  10176 & 4 & 629.3& 1366 & 374.9\\
\hline
 $S_{10}$ & 2498 &   5 & 2273.7 & 34609& 2679.1 \\
\hline
 $S_P$ & 5741 & 4 & 4785.3  & 763650 & 19266.4\\
\hline
\end{tabular}

\medskip

b) $F_4$;

\bigskip

\begin{tabular}{|c|c|c|c|c|c|}
\hline
 Dataset  & Size & Min. length& Avg. length &Max. length & Std. deviation\\
\hline
 $S_1$  & 10165 & 5 & 650.6 & 1388  & 385.4 \\
\hline
 $S_{10}$ & 2566 & 7 & 2791.1 &  28278  & 3234.9 \\
\hline
 $S_P$ & 3821 & 5 & 2430.5 & 160794 & 6491.0\\
\hline
\end{tabular}

\medskip

c) $F_5$;

\end{center}
\caption{Description of the test sets of non-minimal elements  in free
  groups $F_3$, $F_4$ and $F_5$.}
\label{tab:heur_data_non-min}
\end{table}
%

Let ${\mathcal A}$ be one of the  algorithms \DWA,
\HDWA\ or \HPWA.
By an elementary step of the algorithm $\A$,  we mean one
application of a Whitehead automorphism to a given word. Below we
evaluate  the performance of $\A$  with respect to the number of
elementary steps.


Let $N_{s} = N_{s}(\A,S)$ be the average
number of elementary steps  required by $\A$  to find a
minimal element for a  given input $w \in S$, where $S \subset F_n$ is a test set.

By $N_{red} = N_{red}(\A,S)$ we denote the average number of
elementary length-reducing steps required by $\A$ to reduce a
given  element $w \in S$ to a minimal one, so $N_{red}$
is the average number of "productive" steps performed by  $\A$.
 It follows that if $t_1, \ldots, t_l$
are all the length reducing automorphisms found by $\A$ when
executing its routine on an input $w \in S$ then $|w t_1 \ldots
t_l| = |w_{min}|$ and the average value of $l$ is equal to $N_{red}$.

We use values $N_{s}$ and $N_{red}$ as   measures evaluating the performance of the
algorithms.
In addition we record the CPU time $T(w)$ spent by an algorithm to produce a
solution for a particular word $w$.  Since  \HPWA\ is a  probabilistic
algorithm there exists a possibility of producing an incorrect solution. We
measure the error of a probabilistic Whitehead reduction algorithm $\mathcal{A}$
by computing the fraction of elements for which $\mathcal{A}$ failed to return a minimal
element. Let $Sol_{\mathcal{A}}(w) \in F_n$ be a solution produced by algorithm $\mathcal{A}$.
If result is correct, then $|Sol_{\mathcal{A}}(w)| = |w_{min}|$.
The error rate of $\mathcal{A}$ with respect to the test set $S$
\[
E(\mathcal{A}) = \frac{|\{ w \in S \mid |Sol_{\mathcal{A}}(w)| > |w_{min}| \}|}{|S|}.
\]
In all the experiments we have done with the stochastic algorithm \HPWA\ the error rate was zero, i.e.
it has never happened that a non-minimal element has been claimed to be minimal.

First, we experiment with groups of smaller ranks. For elements in  free groups
$F_3$, $F_4$ and $F_5$, algorithm \DWA\ can decide in a practically
acceptable amount of  time on whether an element is minimal or
not. This allows us to obtain the true values of lengths of
minimal elements for each of the input words and access the error rate of probabilistic algorithms.
Results are presented in Tables \ref{tab:HDWA_compare_SP} - \ref{tab:HDWA_compare_S10}, where
\[
T_{avg} = \frac{1}{|S|}\sum_{w \in S}T(w)\\
\]
and $S$ is the corresponding test set.

%
\begin{table}[h]
\begin{center}
\begin{tabular}{|l||r|r||r|r||r|r|}
\hline
 & \multicolumn{2}{|c||}{$N_{s}$} & \multicolumn{2}{|c||}{$N_{red}$} & \multicolumn{2}{|c|}{$T_{avg},s$}  \\
\hline
$\mathcal{A}$ & mean & std  & mean & std & mean& std \\
\hline
\DWA & 360.2 & 267.5& 27.3 & 18.5&  0.11 & 0.46\\
\hline
 \HDWA& 41.2 & 29.1 & 24.2 & 15.0&  0.01 & 0.05 \\
\hline
\HPWA  & 41.1 & 29.6 & 24.2 & 15.0 &  0.01 & 0.05   \\
\hline
\end{tabular}

\medskip

a) $F_3$;

\bigskip

\begin{tabular}{|l||r|r||r|r||r|r|}
\hline
& \multicolumn{2}{|c||}{$N_{s}$} & \multicolumn{2}{|c||}{$N_{red}$} & \multicolumn{2}{|c|}{$T_{avg},s$}  \\
\hline
$\mathcal{A}$ & mean & std  & mean & std & mean& std  \\
\hline
\DWA & 2679.7 & 2356.5& 57.5& 37.3&  2.03& 8.75 \\
\hline
 \HDWA& 118.1 & 114.8& 45.5 & 28.0& 0.08& 0.31  \\
\hline
\HPWA  & 118.3  & 117.1 & 45.5  & 28.0 & 0.07  & 0.26   \\
\hline
\end{tabular}

\medskip

b) $F_4$;

\bigskip

\begin{tabular}{|l||r|r||r|r||r|r|}
\hline
& \multicolumn{2}{|c||}{$N_{s}$} & \multicolumn{2}{|c||}{$N_{red}$} & \multicolumn{2}{|c|}{$T_{avg},s$}  \\
\hline
$\mathcal{A}$ & mean & std  & mean & std & mean& std  \\
 \hline
\DWA & 16319.9 & 20284.53 &79.3 & 52.6& 5.52 &16.4 \\
\hline
 \HDWA& 276.5 & 539.4&58.9 & 38.6 & 0.12& 0.29\\
\hline
\HPWA & 239.2 & 324.8 & 58.0  & 35.6 & 0.08  & 0.16   \\
\hline
\end{tabular}

\medskip

c) $F_5$;

\end{center}
\caption{Comparison of algorithms \DWA, \HDWA\  and \HPWA\ on the test sets of primitive elements $S_P$ in free
  groups $F_3$, $F_4$ and $F_5$, where $N_{s}$ is the average number of elementary steps to find a
minimal element,  $N_{red}$ the average number of length-reducing steps,
 $T_{avg}$ is the average time (in seconds) spent on an input.}
\label{tab:HDWA_compare_SP}
\end{table}
%

%
\begin{table}[h]
\begin{center}
\begin{tabular}{|l||r|r||r|r||r|r|}
\hline
& \multicolumn{2}{|c||}{$N_{s}$} & \multicolumn{2}{|c||}{$N_{red}$} & \multicolumn{2}{|c|}{$T_{avg},s$} \\
\hline
$\mathcal{A}$ & mean & std  & mean & std & mean& std \\
\hline
\DWA  & 129.0 & 33.1 & 2.14 & 1.15 &  0.05 & 0.03\\
\hline
\HDWA & 117.5 & 10.4 & 1.69 & 0.81 & 0.04 & 0.03\\
\hline
\HPWA & 53.6 & 142.4 & 1.70 & 0.82 &  0.011 & 0.01  \\
\hline
\end{tabular}

\medskip

a) $F_3$;

\bigskip

\begin{tabular}{|l||r|r||r|r||r|r|}
\hline
& \multicolumn{2}{|c||}{$N_{s}$} & \multicolumn{2}{|c||}{$N_{red}$} & \multicolumn{2}{|c|}{$T_{avg},s$}\\
\hline
$\mathcal{A}$ & mean & std  & mean & std & mean& std \\
\hline
\DWA &  734.9 & 188.3& 3.24& 1.72& 0.30 & 0.20 \\
\hline
 \HDWA  & 552.1 & 61.1& 2.42 & 1.19&  0.21 & 0.12 \\
\hline
\HPWA  & 140.7  &  387.9 & 2.43 & 1.29& 0.02 & 0.06 \\
\hline
\end{tabular}

\medskip

b) $F_4$;

\bigskip

\begin{tabular}{|l||r|r||r|r||r|r|}
\hline
& \multicolumn{2}{|c||}{$N_{s}$} & \multicolumn{2}{|c||}{$N_{red}$} & \multicolumn{2}{|c|}{$T_{avg},s$} \\
\hline
$\mathcal{A}$ & mean & std  & mean & std & mean& std \\
 \hline
\DWA & 3541.3 & 908.2& 4.28 & 2.19 & 1.45  & 1.05\\
\hline
 \HDWA & 2601.6 & 341.7&  3.29 & 1.73& 0.70& 0.41\\
 \hline
\HPWA &  316.8 & 895.2 & 3.29  & 1.73 &  0.05& 0.06  \\
\hline
\end{tabular}

\medskip

c) $F_5$;

\end{center}
\caption{Comparison of algorithms \DWA, \HDWA\ and \HPWA\ on the test
sets $S_{1}$ in free  groups $F_3$, $F_4$ and $F_5$, where $N_{s}$ is the average number of elementary steps to find a
minimal element,  $N_{red}$ the average number of length-reducing steps,
 $T_{avg}$ is the average time (in seconds) spent on an input.}
\label{tab:HDWA_compare_S1}
\end{table}
%

%
\begin{table}[h]
\begin{center}
\begin{tabular}{|l||r|r||r|r||r|r|}
\hline
& \multicolumn{2}{|c||}{$N_{s}$} & \multicolumn{2}{|c||}{$N_{red}$} & \multicolumn{2}{|c|}{$T_{avg},s$}\\
\hline
$\mathcal{A}$ & mean & std  & mean & std & mean& std \\
\hline
\DWA & 203.5 & 86.5  & 8.45 & 5.49&  0.12 & 0.13 \\
\hline
\HDWA & 124.5  & 13.5& 6.54 & 3.94& 0.05 & 0.03 \\
\hline
\HPWA & 64.7 &  153.9 & 6.54  & 3.94 & 0.02 &  0.01 \\
\hline
\end{tabular}

\medskip

a) $F_3$;

\bigskip

\begin{tabular}{|l||r|r||r|r||r|r|}
\hline
& \multicolumn{2}{|c||}{$N_{s}$} & \multicolumn{2}{|c||}{$N_{red}$} & \multicolumn{2}{|c|}{$T_{avg},s$}\\
\hline
$\mathcal{A}$ & mean & std  & mean & std & mean& std \\
\hline
\DWA  &  1278.6 & 527.5& 17.1 & 10.7&  1.01 & 0.65 \\
\hline
 \HDWA & 569.7 & 67.5 & 11.6 & 7.16& 0.36 & 0.33\\
\hline
\HPWA  &  172.0 &  416.0 &  11.6  & 7.16 & 0.04 & 0.03  \\
\hline
\end{tabular}

\medskip

b) $F_4$;

\bigskip

\begin{tabular}{|l||r|r||r|r||r|r|}
\hline
& \multicolumn{2}{|c||}{$N_{s}$} & \multicolumn{2}{|c||}{$N_{red}$} & \multicolumn{2}{|c|}{$T_{avg},s$}\\
\hline
$\mathcal{A}$ & mean & std  & mean & std & mean& std \\
 \hline
\DWA &  7650.5& 4468.0 & 27.1 &17.8 & 5.87 & 8.81\\
\hline
 \HDWA  & 2650.5 &342.9 & 17.1 & 10.8 & 1.06 & 0.63 \\
 \hline
\HPWA & 360.7 & 904.5 & 16.9 &  10.6 & 0.08 & 0.08 \\
\hline
\end{tabular}

\medskip

c) $F_5$;

\end{center}
\caption{Comparison of algorithms \DWA, \HDWA\ and \HPWA\ on the test
sets $S_{10}$ in free  groups $F_3$, $F_4$ and $F_5$, where $N_{s}$ is the average number of elementary steps to find a
minimal element,  $N_{red}$ the average number of length-reducing steps,
 $T_{avg}$ is the average time (in seconds) spent on an input.}
\label{tab:HDWA_compare_S10}
\end{table}
%
\begin{table}[h]
\begin{center}
\begin{tabular}{|l||r|r||r|r||r|r|}
\hline
& \multicolumn{2}{|c||}{$N_{s}$} & \multicolumn{2}{|c||}{$N_{red}$} & \multicolumn{2}{|c|}{$T_{avg},s$} \\
\hline
 & mean & std  & mean & std & mean& std \\
\hline
\hline
$F_{10}$ & 595.2 & 9195.5& 55.5 & 37.0 & 0.20& 0.58  \\
\hline
\hline
$F_{15}$ & 671.1 & 883.1&  106.1& 55.5 &1.03 & 0.59  \\
\hline
\hline
$F_{20}$ & 736.3 & 874.4&  128.4 & 61.8&  2.80 & 1.41   \\
\hline
\end{tabular}
\end{center}
\caption{Performance of the algorithm \HPWA\ on sets of primitive elements  in free groups $F_{10}$, $F_{15}$ and $F_{20}$,
where $N_{s}$ is the average number of elementary steps to find a
minimal element,  $N_{red}$ the average number of length-reducing steps,
 $T_{avg}$ is the average time (in seconds) spent on an input.}
\label{tab:HPWA_large_ranks}
\end{table}
%

From the tables we can see that both algorithms, \HDWA\ and \HPWA, significantly outperform
\DWA\ on the sets of primitive elements with the error of \HPWA\ being small (actually zero).
This shows that the fast checks are efficient
reduction heuristics. The same picture holds for other sets as well. However, the performance of \HDWA\ deteriorates on
sets $S_1$ and $S_{10}$, where it is much more difficult to decide whether or not an element is minimal.
We  have already mentioned that in the case of
a minimal element all of the Whitehead automorphisms must be applied to confirm that it is indeed minimal. The sizes of the sets of Whitehead elementary
automorphisms in free groups $F_3$, $F_4$ and $F_5$ are $|\Omega_3| = 90$,  $|\Omega_4| = 504$, $|\Omega_5| = 2550$ respectively.
From the tables \ref{tab:HDWA_compare_S1} and \ref{tab:HDWA_compare_S10} we can see that the values of $N_{s}$ in all
cases is just a little greater than the size of $\Omega_n$. This indicates that \HDWA\ spends most of its automorphisms and, therefore,
time, on elements of minimal length. On the contrary, algorithm \HPWA\ seems to be able to avoid exponential blowout by
quickly recognizing minimal elements using the classifier \CLSYS. Note that $N_s$ computed for \HPWA\ is smaller than $|\Omega_n|$
in all experiments.

To show that algorithm \HPWA\ is applicable to groups of large ranks, we perform experiments with primitive elements
in free groups $F_{10}$, $F_{15}$ and $F_{20}$ (see Table \ref{tab:HPWA_large_ranks}). We can see that \HPWA\ was able to find
solutions quickly with $N_s$ growing very slowly with the rank.

\section{Conclusion}

The search heuristics described in \cite{HMM}
can be successfully applied for solving the Whitehead Reduction
problem. Probabilistic algorithm \HPWA\ is very robust and can be
used in groups with large ranks  whereas any other known algorithm
fails to produce similar results due to the fact that the worst case is
inevitable for most inputs. The computational advantage of \HPWA\
increases when the rank of a free group increases. Indeed,  \HPWA\
performs about 11 times faster than \DWA\ in $F_3$ and more
than 60 times faster than in $F_5$.

\section{Acknowledgments}

The authors thank  Alexei G. Miasnikov for his inspiring support of
this work. We also would like to thank  Mike Newman for his valuable
comments and suggestions on this paper.

\bibliographystyle{elsart-harv}
\bibliography{WHybrid}

\begin{thebibliography}{15}
\expandafter\ifx\csname natexlab\endcsname\relax\def\natexlab#1{#1}\fi
\expandafter\ifx\csname url\endcsname\relax
  \def\url#1{\texttt{#1}}\fi
\expandafter\ifx\csname urlprefix\endcsname\relax\def\urlprefix{URL }\fi

\bibitem[{Cohen et~al.(1981)Cohen, Metzler, and Zimmermann}]{Cohen}
Cohen, M., Metzler, W., Zimmermann, A., 1981. What does a basis of $f(a,b)$
  look like? Math. Ann. 257, 435--445.

\bibitem[{Haralick et~al.(2004)Haralick, Miasnikov, and Myasnikov}]{HMM2}
Haralick, R.~M., Miasnikov, A.~D., Myasnikov, A.~G., 2004. Pattern recognition
  approaches to solving combinatorial problems in free groups. Contemporary
  Mathematics 349, 197--213.

\bibitem[{Haralick et~al.(2005)Haralick, Miasnikov, and Myasnikov}]{HMM}
Haralick, R.~M., Miasnikov, A.~D., Myasnikov, A.~G., 2005. {Heuristics for the
  Whitehead Minimization Problem}. {J. Experimental Mathematics} 14~(1), 7--14.

\bibitem[{Kaimanovich et~al.(2005)Kaimanovich, Kapovich, and Schupp}]{KKS}
Kaimanovich, V., Kapovich, I., Schupp, P., 2005. The subadditive ergodic
  theorem and generic stretching factors for free group automorphisms, {Israel
  J. Math.}, to appear, http://arxiv.org/abs/math.GR/0504105.

\bibitem[{Kapovich(2006)}]{K}
Kapovich, I., 2006. Clusters, currents and {Whitehead's} algorithm, preprint,
  http://lanl.arxiv.org/abs/math.GR/0511478.

\bibitem[{Kapovich et~al.(2004)Kapovich, Schupp, and Shpilrain}]{KSS}
Kapovich, I., Schupp, P., Shpilrain, V., 2004. Generic properties of
  {Whitehead's} algorithm and isomorphism rigidity of random one-relator
  groups. Pacific J. Math.To appear.

\bibitem[{Khan(2004)}]{Khan}
Khan, B., 2004. The structure of automorphic conjugacy in the free group of
  rank two. Computational and experimental group theory, Contemp. Math. 349,
  115--196.

\bibitem[{Lee(2003)}]{Lee}
Lee, D., 2003. Counting words of minimum length in an automorphic orbit,
  preprint, http://www.arxiv.org/math.GR/0311410.

\bibitem[{Lyndon and Schupp(1977)}]{LS}
Lyndon, R., Schupp, P., 1977. Combinatorial Group Theory. Series of Modern
  Studies in Math. 89. Springer-Verlag.

\bibitem[{Miasnikov and Shpilrain(2005)}]{MS}
Miasnikov, A., Shpilrain, V., 2005. Automorphic orbits in free groups, {Journal
  of Algebra}, to appear.

\bibitem[{Miasnikov(2004)}]{M}
Miasnikov, A.~D., 2004. Recognition of {Whitehead-minimal} elements in free
  groups of large ranks. Artificial Intelligence and Symbolic Computation
  (Lecture notes in Artificial Intelligence) 3249, 211--221.

\bibitem[{Miasnikov and Myasnikov(2004)}]{MM}
Miasnikov, A.~D., Myasnikov, A.~G., 2004. {Whitehead} method and genetic
  algorithms. Contemporary Mathematics 349, 89--114.

\bibitem[{Moh(1999)}]{moh}
Moh, T.~T., 1999. A public key system with signature and master key functions.
  Communications in Algebra 27~(5), 2207--2222.

\bibitem[{Razborov(1985)}]{Razborov}
Razborov, A., 1985. On systems of equations in a free group. Math. USSR,
  Izvestiya 25~(1), 115--162.

\bibitem[{Whitehead(1936)}]{W}
Whitehead, J.~H.~C., 1936. On equivalent sets of elements in a free group.
  Annals of Mathematic 37, 782 -- 800.

\end{thebibliography}

\end{document}